\documentclass[letterpaper, 10 pt, conference]{ieeeconf}  

\IEEEoverridecommandlockouts                              

\usepackage{graphics} 
\usepackage{epsfig} 
\usepackage{mathptmx} 
\usepackage{times} 
\usepackage{amsmath} 
\usepackage{amssymb}  
\usepackage{cite}               
\usepackage{color}
\usepackage{float}

\newtheorem{definition}{Definition}

\newtheorem{theorem}{Theorem}
\newtheorem{num-example}{Numerical example}
\newtheorem{example}{Example}

\DeclareMathOperator{\sign}{sign}

\newcommand{\omegaE}{\omega_{e}^{\rm{free}}}

\newcommand{\K}{K_{\rm vco}}

\title{\LARGE \bf
	The birth of the global stability theory \\ and the theory of hidden oscillations
}

\author{Kuznetsov~N.V.$^{a}$, 
	Lobachev~M.Y.$^{b}$, 
	Yuldashev~M.V.$^{b}$, 
	Yuldashev~R.V.$^{b}$,\\  
	Kudryashova~E.V.$^{b}$, 
	Kuznetsova~O.A.$^{b}$, 
	Rosenwasser~E.N.$^{c}$, 
	Abramovich~S.M.$^{d}$ 
	\thanks{*This survey is prepared for the invited session 
	``\emph{History of nonlinear systems and control}'' 
	at \emph{the European Control Conference}, Saint Petersburg, 2020. The work is supported by the Leading Scientific Schools of Russia project NSh-2624.2020.1 (sections 1) and the Russian Science Foundation 
	project 19-41-02002 (section 2-3)}
	\thanks{$^{a}$Nikolay V. Kuznetsov is with the Faculty of Mathematics and Mechanics,
		Saint Petersburg State University, Russia, 
		with the Faculty of Mathematical Information Technology,
		University of Jyv\"{a}skyl\"{a}, Finland,
		with the Institute for Problems in Mechanical Engineering RAS, Russia
		{\tt\small nkuznetsov239@gmail.com}}%
	\thanks{$^{b}$Mikhail Y. Lobachev, Marat V. Yuldashev,
	Renat V. Yuldashev, Elena V. Kudryashova, Olga A. Kuznetsova 
	are with the Faculty of Mathematics and Mechanics,
	Saint Petersburg State University, Russia}%
	\thanks{$^{c}$Efim~N.~Rosenwasser is with Saint Petersburg State Marine Technical University, Russia, a participant of the first IFAC World Congress in 1960}
	\thanks{$^{d}$Sergei~M.~Abramovich is with State University of New York at Potsdam, USA}
}

\begin{document}

	\maketitle
	\thispagestyle{empty}
	\pagestyle{empty}

	\begin{abstract}
		
	The first mathematical problems of the global analysis of dynamical models 
	can be traced back to the engineering problem of the Watt governor design.
	Engineering requirements and corresponding mathematical problems led to the fundamental discoveries in the global stability theory.
    Boundaries of global stability in the space of parameters are limited by the birth of oscillations. 
    The excitation of oscillations from unstable equilibria can be easily analysed, while the revealing of oscillations not connected with equilibria is a challenging task being studied in the theory of hidden oscillations.  
	In this survey, a brief history of the first global stability criteria development and corresponding counterexamples with hidden oscillations are discussed.

	\end{abstract}

	\section{INTRODUCTION}
	\label{sec:introduction}

One of the key tasks of the control systems analysis is study of stability and limit dynamic regimes.
A classical example of a mathematical approach to this problem is proposed in I.A.~Vyshnegradsky's work \cite{Vyshnegradsky-1877} on the analysis of governors, 
published in 1877 (see also works of J.C.~Maxwell and A.~Stodola \cite{Maxwell-1868, Stodola-1893, Stodola-1894}).
The design of governors was an important practical task in the XVIII--XIX centuries. 
In 1868 Watt governors were used on about 75 000 steam engines in England alone \cite{Denny-2002}.
However, the absence of a theoretical framework did not allow for controlling the Watt governor's parameters effectively. 
As a result, the operation of steam engines was often unstable, 
and accidents were quite common. 
This problem stimulated the development of stability and control theories.

In his work, Vyshnegradsky, professor of Petersburg Institute of Technology, suggested a mathematical model of system ``steam engine --- Watt governor'' described by the system of ordinary differential equations with one discontinuous nonlinearity. 
For the linearization of this model (obtained by discarding dry friction) 
he determined stability conditions\footnote
{I.A.~Vyshnegradsky's work \cite{Vyshnegradsky-1877} became one of the motivations of A.M.~Lyapunov for further work on rigorous justification of the linearization procedure \cite{Lyapunov-1892} (in 1877 A.M.~Lyapunov was a sophomore in Saint Petersburg University).
}.
Based on engineering reasonings, he conjectured that the obtained conditions are sufficient for the absence of unwanted oscillations and for the transition to sustainable operation with any initial data.

In 1885, H.~L\'{e}aut\'{e} published a paper \cite{Leaute-1885} which showed that governors with dry friction may exhibit nonstationary sustainable regimes.
Later, the famous Russian scientist N.Ye.~Zhukovsky, referring to H.~L\'{e}aut\'{e}'s work, criticized Vyshnegradsky's approach \cite{Zhukovsky-1909} and posed problems of rigorous nonlocal analysis of discontinuous systems and the proof of the Vyshnegradsky's conjecture on the stability of the Watt governor.
This discussion led to the development of the theory of oscillations (studying all posible limit regimes) and the theory of global stability (searching for conditions of the absence of nonstationary limit regimes).

	Significant contribution to the study of oscillations and criteria of its absence 
	was made by A.A.~Andronov's scientific school.
	The monograph ``Theory of oscillations''\cite{AndronovVKh-1937}, first published in 1937, contains the analysis of stability and oscillations of various 
	continuous and discontinuous two-dimensional dynamical models.
		
	Developing this theory further,	A.A.~Andronov and A.G.~Maier studied 
	the three-dimensional nonlinear discontinuous model from \cite{Vyshnegradsky-1877} and proved\footnote{
	The significance of the results was noted when A.A.~Andronov was elected as a full member of Academy of Sciences of the Soviet Union and became the first academician in the field of control theory \cite[p.56]{ANSSSR-1947-eng}.
} that the Vyshnegradsky's conjecture is true \cite{AndronovM-1945, AndronovM-1947,AndronovM-1953}, i.e., that Vyshnegradsky's conditions of local stability imply global stability of the system.

\section{Global stability}
Consider a system of ordinary differential equations
\begin{equation}\label{eq:dot x = f(x)}
\begin{aligned}
\dot{x}=f(x), \quad f:\mathbb{R}^n\to\mathbb{R}^n
\end{aligned}
\end{equation}
and suppose that for any initial state $x_0$ there exists a unique solution
$x(t,\ x_0): x(0,\ x_0)=x_0$, defined on $[0,\ +\infty)$.

\begin{definition}[Global stability]
	\label{def:global stability}
	If any trajectory of system \eqref{eq:dot x = f(x)}	tends to the stationary set,
	then the \emph{system} is called \emph{globally stable}\footnote{
	We use the term ``global stability'' for simplicity of further presentation, while in the literature there are used different terms like 
	``globally asymptotically stable'' \cite[p. 137]{Vidyasagar-1978}, \cite[p. 144]{HaddadC-2011},
	``gradient-like'' \cite[p. 2]{LeonovRS-1992}, \cite[p. 56]{YakubovichLG-2004},
	``quasi-gradient-like'' \cite[p. 2]{LeonovRS-1992}, \cite[p. 56]{YakubovichLG-2004}
	and others,
	reflecting the features of the stationary set and the convergence of trajectories to it.}.
\end{definition}
Note, that the Lyapunov stability of the stationary set in Def.~\ref{def:global stability} is not required. 
An example of a globally stable two-dimensional system with unique unstable equilibrium having a family of homoclinic trajectories can be found in \cite{Vinograd-1957, Hahn-1967}.

Within the framework of global stability study, it is naturally to classify oscillations in control systems as \emph{self-excited} or \emph{hidden} 
\cite{KuznetsovL-2014-IFACWC, BraginVKL-2011, LeonovK-2013-IJBC}.
Basin of attraction of a hidden oscillation in the phase space does not intersect with small neighborhoods of any equilibria, whereas a self-excited oscillation is excited from an unstable equilibrium.
A self-excited oscillation is a \emph{nontrivial} one  
if it does not approach the stationary states
(i.e., $\omega$-limit set of corresponding trajectory does not contain an equilibrium).

The loss of global stability may occur by appearance of either nontrivial 
self-excited oscillation (see, e.g., \cite{Bautin-1949}) or a hidden one.
Self-excited oscillations can be identified by the study of equilibria and computation of trajectories from their vicinities.
However, the revealing of hidden oscillations and obtaining initial data for their computation are challenging problems, which are studied in \emph{the theory of hidden oscillations} 
\cite{Kuznetsov-2018-IFAC,Kuznetsov-2018-MKPU,Kuznetsov-2019-VSPU,Kuznetsov-2019-SCDG,Kuznetsov-2019-UFA,Kuznetsov-2020-TiSU}, 
which represents the genesis of the modern era of Andronov's theory of oscillations.

\subsection{Systems with a single equilibrium}

In 1944, being in Sverdlovsk (now Yekaterinburg), A.I.~Lurie\footnote{
During the wartime in 1941--1944, A.I.~Lurie was in evacuation in Sverdlovsk and chaired the Department of theoretical mechanics at the Ural Industrial Institute \cite{KLurie-2001, Pupyrev-2006}.
The Department regularly held seminars on analytical mechanics and control theory under the guidance of A.I.~Lurie.
} and V.N.~Postnikov published an article \cite{LurieP-1944} with the 
analysis of the global stability of the following model:
\begin{equation}\label{eq:system in the Lurie form}
\begin{aligned}
&\dot{x}=Px+q\varphi(r^Tx),
\end{aligned}
\end{equation}
where $P$ is a matrix, $q$ and $r$ are vectors, and 
$\varphi:\mathbb{R}\to\mathbb{R}$ is a continuous scalar nonlinearity 
such that $\varphi(0)=0$.
Nowadays such models are called {\it Lurie systems} and used to describe various control systems (including the Vyshnegradsky model of the Watt governor). 
In \cite{LurieP-1944} it was suggested to study the global stability of system \eqref{eq:system in the Lurie form} by a Lyapunov function in the form ``quadratic form plus the integral of nonlinearity''.
Later, this class of functions became known as Lyapunov functions of the {\it Lurie-Postnikov form}.
 
The works by Vyshnegradsky, Andronov-Mayer, and Lurie-Postnikov led to the problem of describing a class of Lurie systems for which necessary conditions of stability (i.e., stability of linearized model) coincide with sufficient ones (i.e., global stability of nonlinear model).
In 1949 M.A.~Aizerman, who became acquainted with the Andronov-Mayer results on the stability of the Watt governor at Andronov's seminar in Moscow \cite{Bissel-1998},
formulated the question \cite{Aizerman-1949}: \emph{is the Lurie system with one equilibrium globally stable if the nonlinearity belongs to the sector of linear stability?}
Nowadays, this question is known as the {\it Aizerman's conjecture} on absolute stability.
In 1952, I.G.~Malkin\footnote{
	I.G.~Malkin was the head of the Department of theoretical mechanics of the Ural University (Sverdlovsk).
	He organized a scientific seminar on stability and nonlinear oscillations,
	which was attended by E.A.~Barbashin and N.N.~Krasovsky.
} published an article \cite{Malkin-1952_2} where the method of Lyapunov functions of the Lurie-Postnikov form was developed for the Aizerman's conjecture in the case $n=2$.
In the same year, N.N.~Krasovsky, referring to Malkin's method, presented a counterexample to the Aizerman's conjecture in the case $n=2$ \cite{Krasovsky-1952} with solutions tending to infinity\footnote{See further discussion in \cite{Pliss-1958,LeonovK-2013-IJBC}.}.

Independently, a similar conjecture was later advanced by R.E.~Kalman in 1957,
with the additional requirement that the derivative of nonlinearity
belongs to the linear stability sector~\cite{Kalman-1957}.
Counterexamples with hidden periodic and chaotic oscillations to the Kalman's conjecture can be obtained, for instance, in the four-dimensional Keldysh model of flutter suppression and in model of the Watt governor with a servo motor \cite{Kuznetsov-2018-IFAC,Kuznetsov-2019-SCDG,Kuznetsov-2020-TiSU}\footnote{
Regarding counterexamples with hidden oscillations
 \cite{KuznetsovLS-2011-IFAC,LeonovK-2013-IJBC,KuznetsovL-2014-IFACWC}, R.E.~Kalman wrote that he is, \emph{most certainly, interested in recent developments in the Aizerman 
and (very youthful) Kalman conjecture}.}.


The Lurie's problem and the Aizerman's conjecture stimulated the development 
of general global stability theory.
In 1952, E.A.~Barbashin and N.N.~Krasovsky formulated a general theorem on global stability via Lyapunov functions for autonomous systems of ODEs \cite{BarbashinK-1952}. 
In that paper, the radial unboundedness condition\footnote
{
	$V(x)\to+\infty$ as $||x||\to+\infty$.
} was introduced which allows for conclusions to be made both on local and global stability.

\begin{theorem}(\textit{Barbashin-Krasovsky theorem}\cite{BarbashinK-1952})
	\label{th:BarbashinK}
	Consider system \eqref{eq:dot x = f(x)},
	where $f$ is a continuously differentiable vector-function such that $f(0)=~0$.
	Let $V(x):\mathbb{R}^n\to\mathbb{R}$ be a continuously differentiable function such that:
	
	(i)  $V(x)>0 \quad \forall x\ne0 \quad\text{and} \quad V(0)=0;$

	(ii) $\frac{dV(x(t))}{dt}<0 \quad \forall x\ne0;$
	
	(iii) $V(x)\to+\infty \quad\text{as} \quad||x||\to+\infty.$

	Then any solution tends to the equilibrium $x\equiv0$ 
	(i.e., the system is globally stable)
	and it is Lyapunov stable.
\end{theorem}\medskip

\begin{example}\label{example:example 1} \cite{BarbashinK-1952}
Consider system \eqref{eq:dot x = f(x)} with
\begin{equation}\label{eq:example1}
f(x_1,x_2) = 
     \begin{pmatrix}
	   -\frac{2x_1}{(1+x_1^2)^2}+2x_2\\
	   -\frac{2x_2}{(1+x_1^2)^2}-\frac{2x_1}{(1+x_1^2)^2}
	\end{pmatrix}
\end{equation}
and the Lyapunov function 
 $V(x_1, x_2)=x_2^2+\frac{x_1^2}{1+x_1^2}$.
This function is not radially unbounded, hence, Theorem~\ref{th:BarbashinK} is not applicable for system \eqref{eq:example1} with this Lyapunov function.
In \cite{BarbashinK-1952} it was shown that there is a domain of instability of system \eqref{eq:example1}.

The Lyapunov theorem on asymptotic stability \cite{Lyapunov-1892} provides the existence of transversal level surfaces 
 $\{x\in \mathbb{R}^n: V(x)=c\}$ in the vicinity of the origin, which make trajectories to tend to the origin.  
The additional radial unboundedness condition implies that these level surfaces cover
the whole phase space.

Thus, Example~\ref{example:example 1} demonstrates an importance of the radial unboundedness condition in the global stability analysis.

\end{example}

Various generalizations of this approach were suggested later by J.~LaSalle 
\cite{Lasalle-1960,LaSalleL-1961} and others.
The Barbashin-Krasovsky theorem was modified for non-autonomous systems by V.M.~Matrosov \cite{Matrosov-1962}.

\subsection{Systems with a periodic nonlinearity \& multiple equilibria}
\label{sec:Systems with multiple equilibria}

Consider system \eqref{eq:dot x = f(x)} and suppose that it has a single periodic variable $\sigma$:
\begin{equation}\label{eq:f's periodic property}
\begin{aligned}
&f(z,\sigma+2\pi)=f(z,\sigma) \quad \forall z\in\mathbb{R}^{n-1}, \forall \sigma\in\mathbb{R}.
\end{aligned}
\end{equation}
Then system \eqref{eq:dot x = f(x)} may have multiple equilibria, i.e., if $(z_{\rm eq}, \sigma_{\rm eq})$ is an equilibrium of \eqref{eq:dot x = f(x)} then $(z_{\rm eq}, \sigma_{\rm eq}+~2\pi k)$ is an equilibrium point too for all $k\in\mathbb{Z}$.
System \eqref{eq:dot x = f(x)} can be rewritten in the Lurie form \eqref{eq:system in the Lurie form} with a scalar periodic nonlinearity $\varphi(\sigma)=\varphi(\sigma+2\pi)$.


The first global analysis of such two-dimensional systems with one periodic nonlinearity
was carried out by F.~Tricomi \cite{Tricomi-1933} and A.A.~Andronov \cite{AndronovVKh-1937}.
They used the phase plane analysis method.
In 1959, for such systems Yu.N.~Bakaev suggested to use the Lurie-Postnikov approach
and considered Lyapunov functions of the Lurie-Postnikov form  \cite{Bakaev-1959}: $V(z, \sigma)=~z^THz+\int\limits_0^{\sigma} \varphi(\tau)d\tau$.
It is important to note, that in the cylindrical phase space the Barbashin-Krasovsky theorem cannot be used with Lyapunov functions of the Lurie-Postnikov form for the global stability analysis: the Lyapunov function must be radially unbounded while
$\varphi(\sigma)$ is $2\pi$-periodic function and $V(0, \sigma)=~\int\limits_0^{\sigma} \varphi(\tau) \;d\tau\not\to~+\infty$ as $|\sigma|\to+\infty$ (see also \cite{LeonovKYY-2015-TCAS,Abramovitch-1988}).
Moreover, system \eqref{eq:dot x = f(x)} with nonlinearity \eqref{eq:f's periodic property} may have multiply equilibria.

Later, Bakaev's results were generalized and rigorously justified by G.A.~Leonov 
for system \eqref{eq:dot x = f(x)} 
with the cylindrical phase space \cite{Leonov-PMM-2000, LeonovRS-1992, Leonov-2006-ARC, LeonovK-2014-book}.
\begin{theorem}(\textit{Leonov theorem on global stability for the cylindrical phase space} \cite{Leonov-PMM-2000,LeonovK-2014-book})
	\label{th:Lyapunov-type lemma for PLL}
	Suppose that the stationary set of system \eqref{eq:dot x = f(x)} consists of isolated points, \eqref{eq:f's periodic property} is fulfilled and there exists a continuous function $V(z,\sigma):\mathbb{R}^{n}\to\mathbb{R}$ such that:
	
	(i) $V(z,\sigma+2\pi)=V(z,\sigma) \;\;\;\forall z\in\mathbb{R}^{n-1}, \forall \sigma\in\mathbb{R}$;
	
	(ii) for any solution $x(t)=(z(t), \sigma(t))$ of system \eqref{eq:dot x = f(x)} the function $V(z, \sigma)$ is nonincreasing;
	
	(iii) if $V(z(t), \sigma(t))\equiv V(z(0), \sigma(0))$, then $(z(t),\sigma(t))\equiv (z(0), \sigma(0))$;
	
	(iv)   $V(z,\sigma)+\sigma^2\to+\infty$ as $||z||+|\sigma|\to+\infty$.
	
	Then system \eqref{eq:dot x = f(x)} is globally stable.
\end{theorem}\medskip

\begin{figure*}[t]
	\centering
	\includegraphics[width=0.33\textwidth]{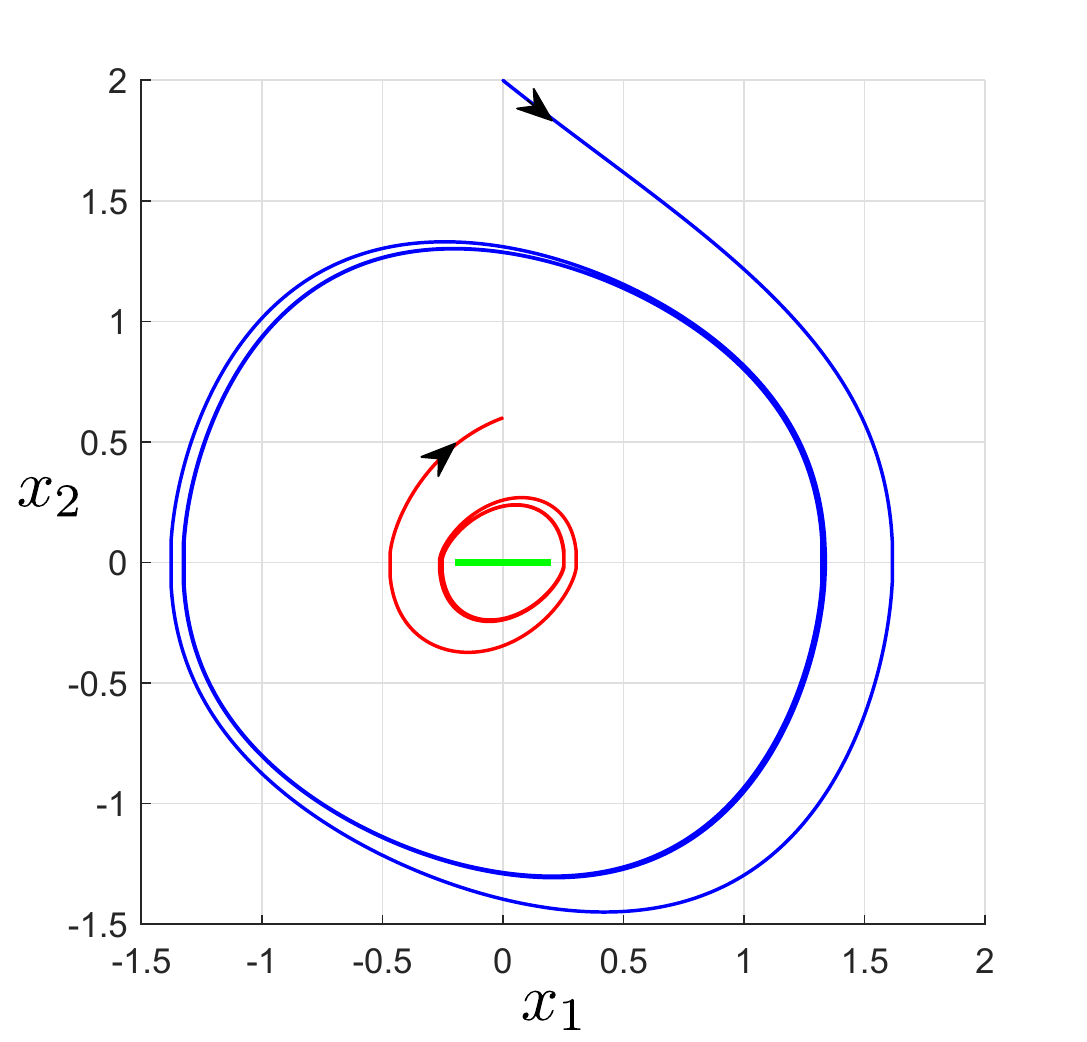}
	\includegraphics[width=0.33\textwidth]{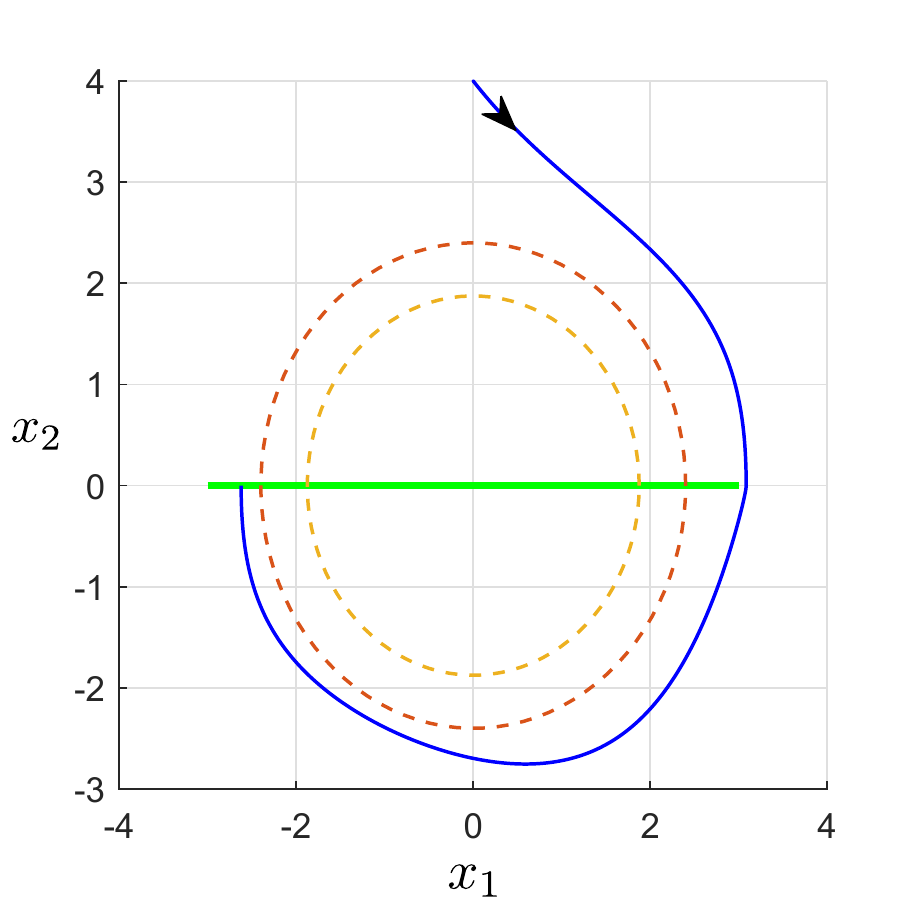}
	\caption{
		Numerical analysis of model \eqref{eq:example3}, $J=k=\kappa=1$.
		Left subfigure: the outer trajectory winds onto the stable limit cycle, 
		the inner trajectory unwinds from the unstable limit cycle and winds onto the stable one (hidden attractor).
		Parameters: $\Phi=0.2$, $\lambda-h\approx-1.2987$ ($\dfrac{\lambda-h}{\sqrt{\Phi\kappa}}\approx-2.9$, i.e., Keldysh's condition \eqref{eq:Keldysh condition} is fulfilled).
		Right subfigure: the outer trajectory approaches the stationary segment, both limit cycles obtained by the harmonic balance method have disappeared (the dash circles).
		Parameters: $\Phi=3$, $\lambda-h\approx-0.937$
		($\dfrac{\lambda-h}{\sqrt{\Phi\kappa}}\approx-2.095$, i.e., Keldysh's condition \eqref{eq:Keldysh condition} is fulfilled).
		}\label{fig}
\end{figure*}

\begin{example} \cite{KuznetsovLYY-2019-DAN,AlexandrovKLNS-2015-IFAC}
	Consider a nonlinear mathematical model of the second-order phase-locked loop 
	with proportionally-integrating filter in the signal's phase space:
	\begin{equation}\label{eq:example2}
	\begin{aligned}
	&\dot z  =  \frac{1}{\tau_1}\sin\sigma, \
	\dot \sigma = \omega_{e}^{\text{free}}
	- \K \left(z + \frac{\tau_2}{\tau_1}\sin\sigma\right),
	\end{aligned}
	\end{equation}
	where $z(t)\in\mathbb{R}$ is a filter state, $\sigma(t)\in\mathbb{R}$ is a phase error, 
	$\K>0$ is a voltage-controlled oscillator (VCO) gain,
	$\omegaE$ is a frequency detuning,
	$\tau_1>0, \tau_2>0$ are parameters.
Here $(\frac{\omegaE}{\K}, 2\pi k), k\in\mathbb{Z}$ are asymptotically stable equilibria and $(\frac{\omegaE}{\K}, \pi+2\pi k), k\in\mathbb{Z}$ are unstable ones.
	Consider the following Lyapunov function of the Lurie-Postnikov form:
	\begin{equation}\label{eq:Lyapunov function from example 2}
	\begin{aligned}
	&V(z, \sigma)= \frac{1}{2}(z-\frac{\omegaE}{\K})^2+
	\frac{1}{\tau_1\K}(1-\cos\sigma).
	\end{aligned}
	\end{equation}
	Its derivative along the trajectories of system \eqref{eq:example2} is
\begin{equation*}
\begin{aligned}
&\dot{V}(z, \sigma)=-
({\tau_2}/{\tau_1^2})\sin^2\sigma<0 \quad\forall \sigma \not\in\{\pi k, \; k\in\mathbb{Z}\}.
\end{aligned}
\end{equation*}

Thus, Lyapunov function \eqref{eq:Lyapunov function from example 2} satisfies the conditions of Theorem~\ref{th:Lyapunov-type lemma for PLL} and system \eqref{eq:example2} is globally stable.
In this case, the model is globally stable for any values of parameters; however, the phase-locked loops with lead-lag filter has only a bounded domain of global stability 
and it is partly determined by the birth of hidden oscillations (see, e.g., \cite{KuznetsovLYY-2017-CNSNS,LeonovK-2013-IJBC,LeonovKYY-2015-TCAS}).
\end{example}

\subsection{Systems with discontinuous nonlinearities}
\label{sec:Systems with discontinuous nonlinearities}

Consider system \eqref{eq:dot x = f(x)}, where $f$ is a piecewise-continuous function with the set of discontinuity points of zero Lebesgue measure.
Discontinuous right-hand side of system \eqref{eq:dot x = f(x)} caused a problem of defining a solution of \eqref{eq:dot x = f(x)} in the discontinuity points.
Thus, it was suggested to consider the solutions as absolutely continuous functions satisfying differential inclusion
\begin{equation}\label{eq:differential inclusion}
\begin{aligned}
\dot{x}\in F(x).
\end{aligned}
\end{equation}
A set $F(x)$ equals to $f(x)$ at continuity points of function $f$.
At discontinuity points $F(x)$ is defined in a special way.
We consider solutions of differential inclusions in Filippov's sense 
\cite{Filippov-1960-IFAC, YakubovichLG-2004}.
 
The generalization of the global stability theorems for the discontinuous systems and the differential inclusions was carried out by A.Kh.~Gelig and G.A.~Leonov \cite{Gelig-1964, Leonov-1971,YakubovichLG-2004} in 1960's-70's\footnote{G.A.~Leonov was a doctoral student at the department chaired by V.A.~Yakubovich and worked on this topic under the guidance of A.Kh.~Gelig 
\cite{AbramovichKL-2015,KuznetsovAFC-2018,AbramovichKN-2018}.}. 
Later, similar results were published in \cite{ShevitzP-1994}.
\begin{theorem}(\textit{Gelig-Leonov theorem on global stability for the differential inclusions} \cite{YakubovichLG-2004, Gelig-1964}).
	\label{th:Lyapunov-type lemma for differential inclusions}
	Let a continuous function $V(x)$ defined in $\mathbb{R}^n$ have the following properties:
	
	(i) $V(x(t))$ is nonincreasing in $t$ for any solution $x(t)$ of \eqref{eq:differential inclusion};
	
	(ii) if the identity $V(x(t))=const$ is valid for all $t\in\mathbb{R}$ and for some solution $x(t)$, bounded when $t\in\mathbb{R}$, then the solution $x(t)$ is a stationary vector;
	
	(iii) $V(x)\to+\infty$ as $||x||\to +\infty$.
	
	Then differential inclusion \eqref{eq:differential inclusion} is globally stable.
\end{theorem}\medskip

Applying ideas of Theorem~\ref{th:Lyapunov-type lemma for differential inclusions} it is also possible to use discontinuous Lyapunov functions for global analysis
(see, e.g., \cite{LeonovKKM-2017} and the proof of the Vyshnegradsky's conjecture on global stability of Watt governor). Also the global stability analysis by discontinuous Lyapunov functions is discussed, e.g., in \cite{PolyakovF-2014, Polyakov-2014}.

\begin{example} \cite{Kuznetsov-2018-IFAC,LeonovK-2018-AIP,LeonovK-2018-DAN}
	Consider system \eqref{eq:system in the Lurie form} with   
$x^T(t)=~(x_1(t),\ x_2(t))\in~\mathbb{R}^2$,
	\begin{equation}\label{eq:example3}
	\begin{aligned}
P=\begin{pmatrix}
	0&J^{-1}\\
	-k&-\mu J^{-1}\\
\end{pmatrix},\
q=\begin{pmatrix}
	0\\
	-1
\end{pmatrix},\
r^T=\begin{pmatrix}
	0&J^{-1}
\end{pmatrix},
	\end{aligned}
	\end{equation}
	representing a two-dimensional (with one degree of freedom) Keldysh model of flutter suppression.
	Here, $J>0$ is the moment of inertia,
	$k>0$ is the stiffness, 
	$\varphi(\sigma)=~(\Phi+~\kappa\sigma^2)~\sign\sigma$ is the nonlinear characteristic of the	hydraulic damper with dry friction, 
	$\Phi$ is the dry friction coefficient,
	$\mu=\lambda-h$, 
	$\ h$ is the proportionality constant,
	$\lambda>0$ and $\kappa>0$ are the damper parameters.
	
M.V.~Keldysh used the harmonic balance method, which is known as an approximate one, to get practically important results on the flatter suppression.
In work \cite{Keldysh-1944} it was stated\footnote{
In his paper \cite{Keldysh-1944} M.V.~Keldysh wrote that \emph{he does not give a rigorous mathematical proof and a number of conclusions are drawn by the intuitive analysis}.
} that under condition
	\begin{equation}\label{eq:Keldysh condition}
	\begin{aligned}
	\lambda-h<-\frac{8}{\pi}\sqrt{2\Phi \kappa/3}\approx - 2.08\sqrt{\Phi \kappa}
	\end{aligned}
	\end{equation}
	system \eqref{eq:example3} has two periodic trajectories (limit cycles), otherwise
	all trajectories of the system converge to the stationary segment.

	The numerical analysis (see Fig.~\ref{fig}) shows that the harmonic balance method can lead to wrong conclusions.
	Two coexisting limit cycles are shown in the left subfigure of Fig.~\ref{fig}. 
	Since there does not exist an open neighbourhood of the stationary segment, which intersects with the outer limit cycle's basin of attraction, then this limit cycle is a hidden oscillation. 
	In the right subfigure of Fig.~\ref{fig} both limit cycles have disappeared and
	the trajectories tend to the stationary segment, while Keldysh's estimate \eqref{eq:Keldysh condition} holds.

	A rigorous study of the Keldysh model \eqref{eq:example3} was performed in \cite{LeonovK-2018-AIP,LeonovK-2018-DAN,KudryashovaKKLM-2019}.
	It was shown that $S=\{x: x_2=0\}$ is a discontinuity manifold, 
	$\Lambda=\{-{\Phi}/{k}\le x_1\le {\Phi}/{k},\ x_2=0\}$ is a stationary
	segment.
	System \eqref{eq:example3} was turned to studying the following differential inclusion:
	\begin{equation*}
	\begin{aligned}
	&\dot{x}\in Px+q\psi(r^Tx),
	\end{aligned} \quad
	\psi(\sigma)=
	\begin{cases}
		&\varphi(\sigma), \quad \sigma\ne0,\\
		&[-\Phi, \Phi], \quad \sigma=0.
	\end{cases}
	\end{equation*}

	Application of the Lyapunov function $V(x_1, x_2) = \frac{1}{2} (kx_1^2+J^{-1}x_2^2)$ and Theorem~\ref{th:Lyapunov-type lemma for differential inclusions} leads to the global stability condition:
	$\lambda-h>-2\sqrt{\Phi \kappa}$.

\end{example}

	\section{CONCLUSIONS}
	\label{sec:conclusion}
	
    As it was noted by Barbashin at the first IFAC World Congress, {\it the methods of constructing Lyapunov functions ... were not sufficiently effective for their use in the investigation of a concrete system}
    \cite{Barbashin-1960-IFAC}.
    To overcome this difficulty, the further development of considered methods has begun and a number of effective global stability criteria was suggested.
    
    One of the first effective criteria of the existence of a Lyapunov function for the systems with smooth right-hand side was obtained in 1954 by Krasovsky \cite{Krasovsky-1954, Krasovsky-1957}.
    His criterion on the existence of a Lyapunov function is based on the stability of the
    symmetrized Jacobi matrix (similar ideas are related to the Markus-Yamabe conjecture \cite{MarkusY-1960} and the corresponding counterexamples with hidden oscillations \cite{Cima-1997, LeonovK-2013-IJBC}).
    Generalizations of the ideas of stability by the first approximation for nonautonomous non-periodic linearizations is a challenging problem because of Perron effects \cite{LeonovK-2007}.

For Lurie systems the development of Lurie-Postnikov approach and the existence of Lyapunov functions are connected with the Popov criterion \cite{Popov-1959, Popov-1960-IFAC, Popov-1961} and famous Kalman-Yakubovich-Popov lemma (KYP lemma) \cite{Yakubovich-1962, Kalman-1963}
(see also \cite{BarabanovGLLMSF-1996}).
In 1959, V.M.~Popov suggested his criterion on absolute stability 
via frequency characteristic (as an electrical engineer he was familiar with frequency characteristics and originally his criterion was not connected with Lyapunov functions).
In 1960, Popov presented this result at the first IFAC World Congress\footnote{
V.M.~Popov's results raised some doubts of M.A.~Aizerman and he asked a young postdoc E.N.~Rosenwasser to find a gap in the Popov's paper.
However, Rosenwasser confirmed the validity of the criterion.
Also, at the first IFAC World Congress, A.I.~Lurie and E.N.~Rosenwasser presented the method of Lyapunov functions construction based on the solvability of so-called Lurie equations \cite{LurieR-1960-IFAC}.} \cite{Popov-1960-IFAC}.

In 1961, Popov proved that the existence of a Lyapunov function of the Lurie-Postnikov form is a sufficient condition for his criterion fulfillment \cite{Popov-1961}.
In 1962, V.A.~Yakubovich formulated the first version of KYP lemma and stated that the converse statement (the necessity of a Lyapunov function existing for the Popov criterion fulfillment) follows from the lemma with some additions \cite{Yakubovich-1962} 
(next year E.N.~Rosenwasser published a paper \cite{Rosenwasser-1963} with detailed explanation of required additions).
Thus, the equivalence of the two approaches was shown.

In 1962, Rosenwasser also noted that the same approach with a Lyapunov function of the quadratic form allows one to obtain the similar criterion for nonautonomous systems \cite{Rosenwasser-1963} and presented the result at Yakubovich's seminar.
Yakubovich generalized this criterion for systems with hysteresis nonlinearities \cite{Yakubovich-1963}\footnote{
	E.N.~Rosenwasser submitted his paper in July, 1962, and V.A.~Yakubovich in September, 1962.
	}.
Nowadays the criterion is known as {\it circle criterion}, and it is sometimes called the Rosenwasser-Yakubovich-Bongiorno criterion \cite{LipkovichF-2015}.

Frequency-domain criteria for systems with discontinuous right-hand sides and for systems with the cylindrical phase space were suggested by Gelig and Leonov (see \cite{YakubovichLG-2004} and refs within).

Remark that the described criteria and its modifications provide only sufficient conditions for global stability, and obtaining necessary and sufficient conditions of global stability is a challenging problem related to the analysis of the boundaries of the global stability in the parameters' space and the birth of oscillations.
While the birth of self-excited oscillations can be effectively identified analytically or numerically, the study of hidden oscillations demands the application of special analytical-numerical methods being developed in \emph{the theory of hidden oscillations}.

	\addtolength{\textheight}{-1cm}   
	

	%
	%
	\section*{ACKNOWLEDGMENT}
	Authors would like to thank F.L.~Chernous'ko, A.L.~Fradkov, 
    R.E.~Kalman~(1930-2016), A.B.~Kurzhanski, V.~R\u{a}svan, and S.N.~Vassilyev
	for valuable comments and discussions.


\end{document}